\documentclass[11pt]{amsart}
\usepackage{longtable, stackrel, multicol, fancybox, stackrel, tcolorbox, mathtools}
\usepackage{bm,amsfonts,amsmath,amssymb,mathrsfs,bbm,amsthm, hyperref} 
\usepackage{graphicx, color, eurosym, eucal,palatino}
\usepackage[margin=1.45in]{geometry}
\usepackage{scrextend}
\changefontsizes{11pt}
\usepackage[normalem]{ulem}

\definecolor{darkred}{rgb}{0.6,0.2,0.2}
\definecolor{darkblue}{rgb}{0.2,0.2,0.6}
\definecolor{darkblue2}{rgb}{0.2,0.2,0.9}
\definecolor{superdarkblue}{rgb}{0.2,0.2,0.3}
\hypersetup{colorlinks=true, linkcolor=darkblue,citecolor=darkblue,
	urlcolor = superdarkblue}
\definecolor{citegreen}{rgb}{0.2,0.2,0.6}

\newcommand\vl[1]{{\color{darkblue} #1}}

\usepackage[english]{babel}
\usepackage{graphicx,amsmath,amssymb,color, enumerate}
\numberwithin{equation}{section}

\usepackage{color}

\newcommand\nb{\nabla}
\newcommand{\beq}{\begin{equation} \begin{split}}
\newcommand{\eeq}{\end{split} \end{equation}}
\newcommand\Sg{\Sigma}
\newcommand\Omg{\Omega}

\makeatletter
\def\section{\@startsection{section}{1}\z@{.9\linespacing\@plus\linespacing}%
	{.7\linespacing} {\fontsize{13}{14}\selectfont\bfseries\centering}}
\def\paragraph{\@startsection{paragraph}{4}%
	\z@{0.3em}{-.5em}%
	{$\bullet$ \ \normalfont\itshape}}

\@addtoreset{equation}{section}
\makeatother

\renewcommand\and{\qquad\text{and}\qquad}

\newcommand\one{\mathbbm{1}}

\newcommand{\comm}[1]{}

\def\sfH{\mathsf{H}}

\def\bm1{\mathbbm{1}}

\def\s{\sigma}
\def\sfh{\mathsf{h}}

\def\p{\partial}

\def\one{\mathbbm{1}}

\def\Im{{\rm Im}\,}

\def\arr{\rightarrow}

\def\tt{\theta}

\def\lm{\lambda}

\def\s{\sigma}

\def\ii{{\mathsf{i}}}
\def\p{\partial}

\def\kp{\kappa}

\def\sfH{\mathsf{H}}
\def\sfh{\mathsf{h}}
\def\one{\mathbbm{1}}

\def\dd{{\mathsf{d}}}

\newcounter{counter_a}
\newenvironment{myenum}{\begin{list}{{\rm(\roman{counter_a})}}%
{\usecounter{counter_a}
\setlength{\itemsep}{1.ex}\setlength{\topsep}{0.8ex}
\setlength{\leftmargin}{5ex}\setlength{\labelwidth}{5ex}}}{\end{list}}

\usepackage[latin1]{inputenc}
\usepackage[T1]{fontenc}

\newcommand{\eg}{{\it e.g.}\,}
\newcommand{\ie}{{\it i.e.}\,}
\newcommand{\cf}{{\it cf.}\,}

\numberwithin{figure}{section}
\numberwithin{equation}{section}
\theoremstyle{plain}
\newtheorem*{thm*}{Theorem}
\newtheorem{thm}{Theorem}[section]

\newtheorem{lem}[thm]{Lemma}
\newtheorem{prop}[thm]{Proposition}

\newtheorem{cor}[thm]{Corollary}

\newtheorem{dfn}[thm]{Definition}
\theoremstyle{remark}
\theoremstyle{plain}

\newcommand{\beu}{\begin{equation*}}
\newcommand{\eeu}{\end{equation*}}
\newcommand{\besu}{\begin{equation*}
\begin{aligned}}
\newcommand{\eesu}{\end{aligned}
\end{equation*}}
\newcommand{\bes}{\begin{equation}
\begin{aligned}}
\newcommand{\ees}{\end{aligned}
\end{equation}}

\newcommand\cA{\mathcal A}
\newcommand\cB{\mathcal B}

\newcommand\cH{\mathcal H}

\renewcommand\frq{\mathfrak q}

\newcommand\ov{\overline}

\newcommand\void[1]{}

\def\ov{\overline}
\def\ran{{\rm ran\,}}

      \def\dC{{\mathbb C}}

      \def\dR{{\mathbb R}}

   \def\dZ{{\mathbb Z}}

   \def\sfH{{\mathsf H}}

\def\cA{{\mathcal A}}   \def\cB{{\mathcal B}}   \def\cC{{\mathcal C}}
      
   \def\cH{{\mathcal H}}   
      
      \def\cO{{\mathcal O}}

\newcommand{\dom}{\mathrm{dom}\,}

\newtheorem{remark}[thm]{Remark}

\newcommand{\Z}{\mathbb{Z}}
\newcommand{\R}{\mathbb{R}}
\newcommand{\C}{\mathbb{C}}
\newcommand{\Ab}{\mathbf{A}}

\newcommand{\eps}{\varepsilon}

\begin{document}

\title[Magnetic isoperimetric inequality]{On the isoperimetric inequality for the magnetic Robin Laplacian with negative boundary parameter}

\author[A. Kachmar]{Ayman Kachmar}
\address[A. Kachmar]{Lebanese University, Department of Mathematics, Nabatiye, Lebanon}
\address{Center for Advanced Mathematical Sciences (CAMS, American University of Beirut)}
\email{akachmar@ul.edu.lb}

\author[V. Lotoreichik]{Vladimir Lotoreichik}
\address[V. Lotoreichik]{Department of Theoretical Physics, Nuclear Physics Institute, 	Czech Academy of Sciences, 25068 \v Re\v z, Czechia}
\email{lotoreichik@ujf.cas.cz}

\date{\today}

\maketitle
\begin{abstract}
We consider the magnetic Robin Laplacian with a negative boundary parameter on a bounded, planar $C^2$-smooth domain. The respective magnetic field is homogeneous. Among a certain class of domains, we prove that the disk maximizes the ground state energy under the fixed perimeter constraint provided that the magnetic field is of moderate strength. This class of domains includes, in particular, all domains that are contained upon translations in the disk of the same perimeter  and all convex centrally symmetric domains.
\end{abstract}

\section{Introduction}\label{sec:int}

Spectral isoperimetric inequalities have a long history in the context of the Laplace operator, dating back to Rayleigh \cite{R} and the celebrated Rayleigh-Faber-Krahn inequality stating that the ball minimizes the Dirichlet ground state energy under the volume constraint \cite{F,K}. 
Ever since,  spectral isoperimetric inequalities  are  the subject of intensive research,  leaving behind  many open questions especially in the presence of magnetic fields.

 Unlike  the case of Neumann boundary condition,  the lowest eigenvalue of  the Laplace operator with a Robin boundary condition does not vanish,  whereby the   inspection of the counterpart of the Rayleigh-Faber-Krahn  inequality makes sense (see \cite{AFK, B, BFNT, D} and the references therein).  Recent contributions show a strong role played by  the  sign of the parameter defining  the Robin  boundary condition.  With  a positive  parameter in hand,   the Robin Laplacian   shares the same kind of isoperimetric inequality with its Dirichlet cousin,  while a negative parameter leads to a  radically   different type of spectral  inequality.  Another central element in the Robin context  is  the dimension of the domain.

The subject of  this paper is on the  challenging question of  the magnetic isoperimetric inequalities.  Our main contribution is a new isoperimetric inequality for the magnetic Laplacian  with  a Robin boundary condition. 

Let  us first  explore some  existing results in  the case without a magnetic field. For a \emph{positive Robin parameter},  the ball minimizes the ground state energy among domains of a fixed volume  \cite{B, D},  very much like the case of  the Dirichlet Laplacian.  However,  for a \emph{negative Robin parameter},  the disk maximizes the  ground state energy among domains of a fixed perimeter~\cite{AFK} and in higher dimensions the ball is known to maximize the ground state energy in the class of convex domains with fixed surface area of the boundary \cite{BFNT}.

Even in the absence of a magnetic field,   the case of  negative Robin parameter is
 mysterious with incomplete  results and 	
a number of unsettled conjectures. For instance, 
it is shown in \cite{FK} for a large negative Robin parameter, the disk is not the maximizer of the ground state energy under fixed area constraint (by providing an example of a non simply connected domain violating the sought property).  It is conjectured in~\cite{AFK} that the disk is still a maximiser among domains of fixed area in the class of simply connected domains.
It is also conjectured in~\cite{AFK} that 
in higher dimensions the ground state energy
is maximized by the ball under fixed surface area of the boundary  without the convexity assumption. 

Apart from the foregoing conjectures,  there are recent interesting results on  the spectral optimization (without a magnetic field) for the lowest Robin eigenvalue in other geometric/topological settings,  like on surfaces and in exterior domains \cite{KL, KrL, KrL2},  and on the  higher Robin eigenvalues as well  \cite{FLa,FLb, GL}.

Isoperimetric  inequalties are rare in the context of the magnetic Laplacian. 
A celebrated result by Erd\H{o}s \cite{E} establishes that {in two dimensions} the {disk} is a minimizer of the magnetic Dirichlet ground state energy under the fixed area constraint provided that the magnetic field is homogeneous. The corresponding question on the {magnetic} Neumann eigenvalue in two dimensions
with  homogeneous magnetic field is still open;  Fournais and Helffer \cite{FH} 
conjecture that the disk is a  maximizer of the ground state energy under fixed area constraint in the class of simply connected domains and support the validity of this conjecture by analysing  asymptotic regimes of weak and strong magnetic fields.  Note that, unlike  the case without magnetic field,  the magnetic ground state energy does not vanish when imposing a Neumann boundary condition,  thereby turning the inspection of the
counterpart of the Rayleigh-Faber-Krahn inequality into  a challenging endeavour. 

Interesting geometric upper bounds on the magnetic Neumann and Dirichlet eigenvalues are derived in~\cite{LS}.   The optimization of the ground state energy for the magnetic Robin Laplacian (including the Neumann case) with the homogeneous magnetic field is still largely  open in the literature.  Intuitively,  the  case with a positive Robin parameter is expected to be effectively  similar to the Dirichlet situation for large boundary parameter
and similar to the Neumann situation for small boundary parameter.

In the present paper we  obtain the two-dimensional isoperimetric inequality
for the magnetic Robin Laplacian with the negative boundary parameter and the homogeneous magnetic field. 
We find that the disk is a \emph{maximizer} within an admissible class of domains with the same perimeter as the {disk} under the assumption that the magnetic field is moderate. Unlike the Neumann setting \cite{FH}, our result is non-asymptotic and holds within a large class of domains. 

Our new spectral inequality is the consequence of  a  tricky  construction of a test function valid in the presence of a {homogeneous} magnetic field with moderate intensity.    This test function depends on the distance to the boundary only and in its construction we rely on the fact that for the weak magnetic field the ground-state eigenfunction of the magnetic Robin Laplacian on the disk with a negative boundary parameter is radial and the lowest eigenvalue on the disk is negative. The class of admissible domains
is characterized by a purely geometric condition somehow related to the classical  optimization of the moment of inertia  of curves  \cite{H02}. This class includes all domains that are contained upon translations in the disk of the same perimeter and  all convex centrally
symmetric domains.  The class sounds rather generic but  it remains an open question whether there are simply connected domains outside it.

The body of the paper consists of four sections and three appendices.
Section~\ref{sec:L-uf} introduces the Robin Laplacian we are concerned with.
The case of the disk is analysed in Section~\ref{sec:d}.  
Our main result on the isoperimetric inequality, Theorem~\ref{thm:ipi}, is contained in Section~\ref{sec:ipi}. 
In Section~\ref{sec:sa} we discuss the isoperimetric inequality in the context of large coupling asymptotics for general domains with smooth boundaries. 
In Appendices~\ref{app:clsd_smb} and \ref{sec:cont}, we collect standard arguments related to the definition of the Robin Laplacian and the continuity of its eigenvalues. Finally, a standard result on the magnetic Neumann Laplacian is recalled in Appendix~\ref{app:C}.
\section{The Robin Laplacian with a homogeneous magnetic field}\label{sec:L-uf}

Consider a bounded simply connected planar domain $\Omega\subset\R^2$ with a $C^2$-smooth boundary $\partial\Omega$ having the length
\begin{equation}\label{eq:def-L}
|\partial\Omega|=L\,.
\end{equation}
Given two parameters $b\geq 0$ (the intensity of the magnetic field) and $\beta\le0$ (the Robin parameter), consider the closed, densely defined symmetric and semi-bounded quadratic form
\begin{equation}\label{eq:qf}
	\frq_{\Omega}^{\beta,b}[u] :=\|(\nabla-\ii b\Ab)u\|^2_{L^2(\Omg;\dC^2)}+\beta\|u|_{\p\Omg}\|_{L^2(\partial\Omega)}^2,\quad \dom\frq_{\Omg}^{\beta,b} := H^1(\Omg),
\end{equation}
where the vector potential $\Ab$ is defined by
\begin{equation}\label{eq:mp}
\Ab(x):=\frac12(-x_2,x_1),\quad\big( x=(x_1,x_2)\big)\,.
\end{equation}
For the convenience of the reader we provide in Appendix~\ref{app:clsd_smb} a proof of closedness and semi-boundedness of the form $\frq_\Omg^{\beta,b}$.
\begin{dfn}\label{def:Op}
	The magnetic Robin Laplacian $\sfH_{\Omg}^{\beta,b}$ in the Hilbert space $L^2(\Omg)$ is defined as the unique self-adjoint operator associated with the quadratic form $\frq_{\Omg}^{\beta,b}$ via the first representation theorem~\cite[Thm. VI.2.1]{Kato}.
\end{dfn}
The operator $\sfH_{\Omg}^{\beta,b}$ is characterised by
\[
\begin{split}
\dom\sfH_{\Omg}^{\beta,b} &\! = \!
\big\{
u\!\in\! H^1(\Omg)\colon
\exists\, w\in L^2(\Omg):
\frq_{\Omg}^{\beta,b}[u,v] = (w, 
	v)_{L^2(\Omg)}, \forall\,v\in\dom\frq_{\Omg}^{\beta,b}
\big\},\\
\sfH_{\Omg}^{\beta,b} u&=  w;
\end{split}
\]
here the function $w$ in the characterisation of the operator domain is unique if it exists and hence the operator $\sfH^{\beta,b}_\Omg$ is well defined.
We get then integrating by parts that
\begin{equation}\label{eq:operator}
\begin{split}
\dom\sfH_{\Omg}^{\beta,b} & = \big\{u\in H^1(\Omg)\colon (\nabla-\ii b\Ab)^2u\in L^2(\Omg),\,\, 
\nu\cdot(\nabla-\ii b\Ab)u=\beta u\,\,{\rm on}\,\partial\Omega\big\},\\
\sfH_{\Omg}^{\beta,b}u& =-(\nabla-\ii b\Ab)^2u=-\Delta u+2\ii b\Ab\cdot\nabla u+b^2|\Ab|^2u\,,
\end{split}
\end{equation}
where $\nu$ is the unit inward normal vector of $\partial\Omega$. For all $\beta < 0$, taking into account the smoothness of the boundary
and that $\nu\cdot(\nabla -\ii b\Ab)u|_{\p\Omg} \in H^{1/2}(\p\Omg)$ for all $u\in\dom\sfH_\Omg^{\beta,b}$,   
the elliptic regularity estimates (\cf \cite[Thm. 4.18\,(ii)]{McL})
 yield that $\dom\sfH_{\Omg}^{\beta,b}$
consists of functions in the Sobolev space $H^2(\Omega)$ that satisfy the (magnetic) Robin condition $\nu\cdot(\nabla-ib\Ab)u=\beta u$ on $\partial\Omega$.

It follows from the compact embedding of $H^1(\Omg)$ into $L^2(\Omg)$ that the spectrum of $\sfH_{\Omg}^{\beta,b}$ is purely discrete.
The lowest eigenvalue of the self-adjoint operator $\sfH_\Omg^{\beta,b}$ is characterised by the min-max principle
\begin{equation}\label{eq:ev}
	\lambda_1^{\beta,b}(\Omega) {:=} \inf_{u\in H^1(\Omega)\setminus\{0\}} \frac{ {\frq_{\Omega}^{\beta,b}[u]}}{\|u\|^2_{L^2(\Omg)}}\,.
\end{equation}  
Since $\Omega$ is  simply connected, the eigenvalue $\lambda_1^{\beta,b}(\Omg)$ is independent of the choice of the  vector potential $\Ab$ of the magnetic field. This is a consequence of  invariance under  gauge transformations; if $\Ab'\in H^1(\Omega;\R^2)$ and ${\rm curl}\Ab'=1$, then $\Ab'=\Ab+\nabla\phi$ for a function 
$\phi\in H^2(\Omg)$ (\cf \cite[Props D.1.1 and D.2.1]{FH-b}),  and in turn
\[ \frac{\| (\nabla-\ii b\Ab')u\|^2_{L^2(\Omg;\C^2)}+\beta\|u|_{\partial\Omg}\|_{L^2(\partial\Omg)}^2}{\|u\|_{L^2(\Omg)}^2} =\frac{\frq_{\Omg}^{\beta,b}[e^{-\ii b\phi} u]}{\|e^{-\ii b\phi}u\|_{L^2(\Omg)}^2}\,.\]

Since the quadratic form $\frq_{\Omega}^{\beta,b}[u]$ is continuous with respect to $(\beta,b)$ uniformly in $u$,  a classical theorem yields that the  eigenvalue $\lambda_1^{\beta,b}(\Omega)$ depends  continuously on $(\beta,b)\in\ov{\dR_-}\times\ov{\dR_+}$, where $\dR_- = (-\infty,0)$ and $\dR_+ = (0,\infty)$.  For  convenience, we give a short reminder of this standard material in Appendix~\ref{sec:cont}. 
We introduce the following constant
\begin{equation}\label{eq:beta-c}
\beta_c(b,\Omega):=\sup\{ \beta\in\dR_-\colon\lambda_1^{\beta,b}(\Omega)<0\} < 0\,.
\end{equation}
It should be mentioned that in view of~\eqref{eq:ev} applied to the constant test function one can easily check that $\lm_1^{\beta,b}(\Omg)$ is indeed negative for $\beta < 0$ large by absolute value.
Notice that
\begin{equation}\label{eq:beta-c*}
\beta_c(0,\Omega)=0\quad{~\rm and}\quad\beta_c(b,\Omega)<0~{\rm for~}b>0\,,
\end{equation}
since $\lim\limits_{\beta\to0^-}\lambda_1^{\beta,b}(\Omega)=\lambda_1^{0,b}(\Omega)$. In fact, $\lambda_1^{0,b}(\Omega)$ is the magnetic Neumann eigenvalue; it is non-negative and vanishes if, and only if, $b=0$ (see Appendix~\ref{app:C}).

Since $\lambda_1^{\beta,b}(\Omega)$ is a monotone function of $\beta$, we observe that
\begin{equation}\label{eq:cond:lambda<0}
\lambda_1^{\beta,b}(\Omg)<0\quad{\text{if, and only if,}} \quad \beta<\beta_{\rm c}(b,\Omega)\,.
\end{equation}

\section{The case of the {disk}}\label{sec:d}
In this section we analyse the magnetic Robin Laplacian with a negative boundary parameter on the disk. Related
analysis of the magnetic Laplacian on the disk appears in the literature for the Dirichlet (see \eg~\cite{Son}) and the Neumann (see \eg~\cite{FS}) boundary conditions. 	
Consider a fixed constant $R>0$ and the disk
\begin{equation}\label{eq:disc}
\cB= \cB_R :=\{x\in\R^2\colon|x|<R\}\,.
\end{equation}
We can express the $L^2$-norm in $L^2(\cB)$ and the quadratic form 
{$\frq_{\cB}^{\beta,b}$ in polar coordinates,
\begin{equation}
\label{eq:qf-polar}
\begin{aligned}
\|u\|^2_{L^2(\Omg)}&\!=\!\int_0^{2\pi}\int_0^R |u|^2r\dd r\dd\theta, \\
\frq_{\cB}^{\beta,b}[u]&\!=\!\int_0^{2\pi}\int_0^R \left(|\p_r u|^2\!+\!\frac1{r^2}\Big|\p_\tt u-\frac{\ii br^2}2u\Big|^2\right)r\dd r\dd\tt\!+\!\beta R \int_0^{2\pi}|u(R,\tt) |^2\dd\theta,
\end{aligned}
\end{equation}
where in order to represent $\frq_{\cB}^{\beta,b}$ we used the expression for the magnetic gradient
\[
\nabla -\ii b{\bf A} ={\bf e_r} \p_r  + {\bf e}_\tt\left(\frac{\p_\tt}{r} -\frac{\ii b r}{2}\right),
\]
in which
the moving frame $({\bf e}_r,{\bf e}_\tt)$ 
associated with the polar coordinates is defined by
\[
	{\bf e}_r := \begin{pmatrix}\cos\tt\\
	\sin\tt\end{pmatrix},\qquad {\bf e}_\tt := \begin{pmatrix}-\sin\tt\\\cos\tt\end{pmatrix}.
\]
\subsection*{Fiber operators}
We can separate variables by working in polar coordinates and doing the  Fourier transform with respect to the angular variable.
To this aim we consider the complete family of mutually orthogonal projections in the Hilbert space $L^2(\cB)$
\[	
	(\Pi_m u)(r,\tt) = \frac{1}{2\pi}e^{\ii m \tt}\int_0^{2\pi} u(r,\tt') e^{-\ii m\tt'}\dd\tt',\qquad m\in\dZ.
\]
Upon natural identification of $\ran\Pi_m$ and $L^2((0,R);r\dd r)$, this family of projections induces the orthogonal decomposition
\begin{equation}\label{eq:orthospace}
	L^2(\cB) \simeq \bigoplus_{m\in\dZ} L^2((0,R);r\dd r).
\end{equation}
Using the representation~\eqref{eq:qf-polar} of the quadratic form
$\frq_{\cB}^{\beta,b}$ in polar coordinates
we arrive at the family of closed, densely defined, symmetric and semi-bounded quadratic forms ($m\in\dZ$) in the Hilbert space $L^2((0,R);r\dd r)$
\begin{equation}\label{eq:qf-1d}
\begin{aligned}
\frq_{m,R}^{\beta,b}[f] &:= \frq^{\beta,b}_{\cB}\left[\frac{f(r)e^{\ii m\tt}}{\sqrt{2\pi}}\right] \\
&=\int_0^R \left(|f'(r)|^2+\frac1{r^2}\left(m-\frac{br^2}{2} \right)^2|f|^2 \right)r\dd r+\beta R|f(R)|^2,\\
\dom\frq_{m,R}^{\beta,b} &:=
\left\{f\in L^2((0,R);r\dd r)\colon f(r)e^{\ii m \tt}\in H^1(\cB)\right\}\\
&= \big\{f\colon f,f',mr^{-1}f \in L^2((0,R);r\dd r)\big\}.
\end{aligned}
\end{equation}
Employing the characterisation of the operator $\sfH_\cB^{\beta,b}$ in~\eqref{eq:operator} one can easily check that
\[
	\Pi_m(\dom\sfH_{\cB}^{\beta,b})\subset\dom\sfH_{\cB}^{\beta,b}\quad\text{and}\quad\sfH^{\beta,b}_{\cB}\big(\ran\Pi_m\cap\dom\sfH^{\beta,b}_\cB))\subset\Pi_m(L^2(\cB)).
\]
Let $\sfH^{\beta,b}_{m,R}$ be the self-adjoint fiber operator in the Hilbert space $L^2((0,R);r\dd r)$ associated with the form $\frq_{m,R}^{\beta,b}$, via the first representation theorem. 

\begin{remark}	
The aim of this remark is to characterise the
fiber operators $\sfH_{m,R}^{\beta,b}$. This characterisation essentially follows from the analysis of the Bessel-type operators on an interval; see \eg \cite{AB15, BG, GPS, KL10}.
To this aim we associate with the differential expression
\[
	\ell_m := -\frac{\dd^2}{\dd r^2}-\frac1r\frac{\dd}{\dd r}+\frac1{r^2}\left(m-\frac{br^2}{2} \right)^2,\qquad m\in\dZ\,,
\] 	
the self-adjoint Sturm-Liouville operator
\[
\begin{aligned}
	\cH_{m,R}^{\beta,b} f & :=\ell_m f,\\
	\dom\cH_{m,R}^{\beta,b}& := \Big\{f\colon f,\ell_m f\in L^2((0,R);r\dd r)\\
	&\qquad\qquad\qquad f'(R) = -\beta f(R)~\text{and}~\lim_{r\arr0^+} \tfrac{f(r)}{\ln r} = 0~\text{for}~m= 0\Big\},
\end{aligned}	
\]
acting in the Hilbert space $L^2((0,R);r\dd r)$.
Using the expansions of the type~\cite[Thm. 2.2]{KL10} one can check the inclusion~$\dom\cH^{\beta,b}_{m,R}\subset\dom\frq^{\beta,b}_{m,R}$.
Integrating by parts for any $f\in\dom\cH^{\beta,b}_{m,R}\subset\dom\frq^{\beta,b}_{m,R}$ and $\phi \in \dom\frq^{\beta,b}_{m,R}$ we observe that
\[
\begin{aligned}
\frq_{m,R}^{\beta,b}[f,\phi]&=\int_0^R   \Big(-\frac1r\big(rf'(r)\big)'  +\frac1{r^2}\Big(m-\frac{br^2}{2} \Big)^2 f(r) \Big) \overline{\phi(r)}\,r\dd r\\
&\qquad\qquad-\lim_{r\arr0^+} rf'(r)\ov{\phi(r)}+Rf'(R)\overline{\phi(R)}+\beta Rf(R)\ov{\phi(R)}\\
&= \int_0^R(\cH^{\beta,b}_{m,R}f)(r)\ov{\phi(r)}r\dd r,
\end{aligned}\]
where
	$\lim\limits_{r\arr0^+} rf'(r)\ov{\phi(r)} = 0$ thanks to combination of the expansions~\cite[Thm. 2.2]{KL10} and of~\cite[Eq. (4.14)]{GPS} adapted to our setting, see also~\cite[Prop. 3.2\,(i)]{AB15}.
Hence, the first representation theorem yields that $\cH^{\beta,b}_{m,R}\subset\sfH^{\beta,b}_{m,R}$
and since both operators are self-adjoint, they coincide. 
\end{remark}

In view of the identification between the spaces $\ran\Pi_m$ and $L^2((0,R);r\dd r)$ it follows from the above construction that $\sfH^{\beta,b}_{m,R}$ can be identified with $\sfH_\cB^{\beta,b}|_{\Pi_m(\dom\sfH_\cB^{\beta,b})}$ on $\Pi_m(L^2(\cB))$.
Hence, according to~\cite[\S 1.4]{S12}
we end up with the orthogonal decomposition
\begin{equation}\label{eq:ortho}
	\sfH_{\cB}^{\beta,b} \simeq\bigoplus_{m\in\dZ} \sfH^{\beta,b}_{m,R}
\end{equation}
with respect to ~\eqref{eq:orthospace}.
From the above decomposition and the fact that the spectrum of $\sfH_{\cB}^{\beta,b}$ is purely discrete it follows that the spectra of the fiber operators are also purely discrete.
The lowest eigenvalues of the fiber operators are characterised by
\begin{equation}\label{eq:ev-1d}
	\mu_{1,m}^{\beta,b}(R)=\inf_{f\in\dom\frq_{m,R}^{\beta,b}\setminus\{0\}}
		\frac{ \frq_{m,R}^{\beta,b}[f]}{\int_0^R|f|^2r\dd r}.
\end{equation} 
Moreover, if a function $f\in\dom\frq_{m,R}^{\beta,b}$ minimizes the Rayleigh quotient in \eqref{eq:ev-1d}, then it is an eigenfunction associated with eigenvalue $\mu_{1,m}^{\beta,b}(R)$ (see~\cite[\S.~10.2, Thm. 1]{BS87}).
Relying on the orthogonal decomposition~\eqref{eq:ortho}
the lowest eigenvalue of the magnetic Robin Laplacian $\sfH_{\cB}^{\beta,b}$ is given by
\begin{equation}\label{eq:ev-disc-p*}
\lambda_1^{\beta,b}(\cB)=\inf_{m\in\Z}\mu_{1,m}^{\beta,b}(R)\,.
\end{equation}
In the next proposition we use Sturm-Liouville theory to show that the eigenvalues $\big(\mu_{1,m}^{\beta,b}(R)\big)_{m\in\Z}$ are  all simple. This claim is analogous to
\cite[Lem. 2.2]{BPT98}, where only the Neumann boundary condition is covered.
\begin{prop}\label{prop:StLi}
For all $m\in\Z$, the lowest eigenvalue $\mu_{1,m}^{\beta,b}(R)$ of $\sfH_{m,R}^{\beta,b}$ is simple and
and the respective normalized eigenfunction $f_m$ can be chosen positive on $(0,R)$.
\end{prop}
\begin{proof}
Pick a normalized ground state $u_m$ of $\sfH^{\beta,b}_{m,R}$. It is easy to see that this ground-state can be chosen to be real-valued. Let $f_m=|u_m|$, then $f_m$ is a normalized ground state too,  since
\[ \int_0^R |f_m|^2r\dd r=\int_0^R|u_m|^2r\dd r=1 \quad{\rm and}\quad \frq_{m,R}^{\beta,b}[f_m] = \frq_{m,R}^{\beta,b}[u_m]= \mu_{1,m}^{\beta,b}(R)\,. \]
In particular, we have $f_m\in\dom\sfH^{\beta,b}_{m,R}$ and hence 
	$f_m$ is continuously differentiable on $(0,R)$.
If  $f_m$ vanishes at some point $r_0\in(0,R)$, then $f_m'(r_0)=0$ because $f_m\geq0$, hence
\[ 
\begin{cases}
\sfH^{\beta,b}_{m,R}f_m=-f''_m-\frac1rf'_m+\frac1{r^2}\left(m-\frac{br^2}2\right)^2f_m=\mu_{1,m}^{\beta,b}(R) f_m~{\rm on}~(0,R)\\
f_m(r_0)=f'_m(r_0)=0\quad{\rm and}\quad f_m'(R)=-\beta f_m(R)\end{cases}\,,\]
which yields  $f_m=0$ on $[r_0,R]$, by Cauchy's uniqueness theorem for ODE.  The same argument yields $f_m=0$ on $(0,r_0]$, hence $f_m\equiv0$ which is impossible. Therefore, we must have $f_m>0$ everywhere on $(0,R)$ and hence $u_m$
is strictly sign definite on $(0,R)$. Consequently,  it is impossible to find two orthogonal eigenfunctions corresponding to $\mu_{1,m}^{\beta,b}(R)$.
\end{proof}

\subsection*{Structure of {the ground state}}
It follows from the orthogonal decomposition~\eqref{eq:ortho} that
if $m_\star\in\dZ$ is such that 
\[
	\lambda_1^{\beta,b}(\cB)=\mu_{1,m_\star}^{\beta,b}(R)\,,
\]
then an eigenfunction represented by
\begin{equation}\label{eq:ground-state}
	u_1^{\beta,b}(r,\tt) = f_\star(r)e^{\ii m_\star\tt},
\end{equation}
with $f_\star:=f_{m_\star}$ being the positive normalized ground
state  of $\sfH_{m_\star,R}^{\beta,b}$\,, is associated to the lowest eigenvalue $\lm_1^{\beta,b}(\cB)$ of the operator $\sfH_{\cB}^{\beta,b}$. 

\begin{prop}\label{prop:ev-1d}
Let the self-adjoint operator $\sfH_{\cB}^{\beta,b}$ be associated with the quadratic form $\frq_{\cB}^{\beta,b}$ in~\eqref{eq:qf} as in Definition~\ref{def:Op}.
Then the following hold.
\begin{myenum}
\item There exist $m_\star=m_\star(\beta,b,R)\in\Z$ such that $|m_\star(\beta,b,R)|\le b R^2$ and
\[
	\lambda_{1}^{\beta,b}(\cB) =  \mu_{1,m_\star}^{\beta,b}(R).
\]
\item If $bR^2 < 1$, then to the lowest eigenvalue 
$\lm_1^{\beta,b}(\cB)$ of $\sfH_{\cB}^{\beta,b}$ 
corresponds a radial eigenfunction.
\end{myenum}
\end{prop}
\begin{proof}
(i) 
Suppose that $|m| > bR^2$.
Let us introduce the potential  $V_{m}^{b,R}(r)$\break$=\frac1{r^2}\Big(m-\frac{br^2}{2} \Big)^2$, $r\in (0,R)$. Notice that
\begin{equation}\label{eq:potential}
\begin{aligned}
	V_{m}^{b,R}(r) 
	&= \frac{b^2r^2}{4} + \frac{m^2}{r^2}-mb \\
	&> \frac{b^2r^2}{4} + \frac{|m|bR^2}{r^2}-mb\\
	&\ge
	\frac{b^2r^2}{4} + |m|b-mb \ge V_{0}^{b,R}(r),
\end{aligned}
\end{equation}
where we used that $|m|> bR^2$ in the second step. It follows from~\eqref{eq:qf-1d} that   $\dom\frq_{0,R}^{\beta,b}\supseteq\dom\frq_{m,R}^{\beta,b}$ and thanks to~\eqref{eq:potential} we have $\frq_{0,R}^{\beta,b}[f] < \frq_{m,R}^{\beta,b}[f]$ for all $f\in\dom\frq_{m,R}^{\beta,b}$. Hence,  the characterisation~\eqref{eq:ev-1d} implies
\[
	\mu_{1,m}^{\beta,b}(R) > \mu_{1,0}^{\beta,b}(R).
\]
The claim follows from the above inequality combined with~\eqref{eq:ev-disc-p*}.

\noindent (ii) It follows from (i) that $bR^2 < 1$ implies $\lm_1^{\beta,b}(\cB) = \mu_{1,0}^{\beta,b}(R)$. Hence,~\eqref{eq:ground-state} yields that a radial eigenfunction corresponds to the lowest eigenvalue $\lm_1^{\beta,b}(\cB)$ of the operator $\sfH_{\cB}^{\beta,b}$. 
\end{proof}

\begin{remark}\label{rem:b*-rad}
Let us introduce the following set
\begin{equation}\label{eq:b*-rad}
\cA =\{(\beta,b)\in\R_-\times\R_+\colon\lambda_1^{\beta,b}(\cB)<0 {\rm ~and~} \sfH_{\cB}^{\beta,b}\text{ has a radial ground state}\}\,.
\end{equation}
By \eqref{eq:cond:lambda<0} and Proposition~\ref{prop:ev-1d}, $\cA\not=\emptyset$; in fact, if $bR^2<1$ and $\beta<\beta_{\rm c}(b,\cB)$, then $(\beta,b)\in\cA$.
\end{remark}

In the case where a radial ground state exists, we recall further regularity properties that will be used in our proof of the isoperimetric inequality.

\begin{prop}\label{prop:ef-rad}
Assume that $(\beta,b)\in\cA\subset\dR_-\times\dR_+$ where the set $\cA$ is as in~\eqref{eq:b*-rad}. Let
$u_1^{\beta,b}(x) = f_\star(|x|)$ be the radial ground-state
of the operator $\sfH_{\cB}^{\beta,b}$ 
corresponding to its lowest eigenvalue $\lambda_1^{\beta,b}(\cB)< 0$ represented as in~\eqref{eq:ground-state}. Then $f_\star\in C^\infty([0,R])$, $f_\star'(0)=0$ and $f_\star > 0$ on $(0,R)$.
\end{prop}

\begin{proof}
By the elliptic estimates~\cite[Thm. 4.18\,(ii)]{McL}, $u_1^{\beta,b}\in C^\infty(\ov{\cB})$. 
Clearly,  $f_\star\in C^\infty([0,R])$ since
\[ f_\star(r)=u_1^{\beta,b}(r,0),\qquad r\in[0,R].\]
Furthermore, $f_\star'(0) = 0$ because $f_\star(r)= u_1^{\beta,b}(r,0)=u_1^{\beta,b}(-r,0)$, for all $r\in [0,R]$. Finally,
it follows from Proposition~\ref{prop:StLi} and the representation~\eqref{eq:ground-state} with $m_\star = 0$ that $f_\star(r) > 0$ for all $r\in(0,R)$.
\end{proof}
\begin{remark}
	With additional efforts one can show that $f_\star(0), f_\star(R) > 0$ in the above proposition, but this is not needed for our analysis.
\end{remark}

\subsection*{Estimate of $\beta_{\rm c}(b,\cB)$}
In the next proposition we use the constant test function in order to estimate the critical boundary parameter $\beta_{\rm c}(b,\cB)$.
\begin{prop}\label{prop:estimate_beta}
	Let $b > 0$ be arbitrary. Then the critical boundary parameter $\beta_{\rm c}(b,\cB)$ defined as in~\eqref{eq:beta-c} satisfies
	\[
	\beta_{\rm c}(b,\cB) \ge -\frac{R^3b^2}{16}.
	\]
\end{prop}
\begin{proof}
	Substituting the characteristic function $\one_\cB$ of the disk $\cB$ into the min-max principle~\eqref{eq:ev}  we find that
	\[
		\lm_1^{\beta,b}(\cB) \le \frac{\frq_{\cB}^{\beta,b}[\one_\cB]}{\|\one_\cB\|^2_{L^2(\cB)}} 
		=
		\frac{\frac{\pi b^2}{2} \int_0^R r^3\dd r + \beta|\p\cB|}{|\cB|}
		=
		 \frac{\frac{ R^3}{8}b^2 +2  \beta}{R}. 
	\]
	Hence, for all $\beta < -\frac{R^3}{16}b^2$ we have $\lm_1^{\beta,b}(\cB) < 0$ and the claim follows. 
\end{proof}
\section{An isoperimetric inequality}\label{sec:ipi}
In this section we formulate and prove an isoperimetric inequality for the lowest eigenvalue of the magnetic Robin Laplacian with a negative boundary parameter. The argument is inspired by the proof of a similar inequality for the non-magnetic Robin Laplacian~\cite[Thm. 2]{AFK} and relies on the method of parallel coordinates. In order to include the magnetic term into consideration an additional geometric assumption will be imposed. 
 
Let $\cB\subset\dR^2$ be a disk of the same perimeter $L > 0$ as a $C^2$-smooth simply connected domain $\Omg\subset\dR^2$. We denote by $R=\frac{L}{2\pi} >0$ the radius of $\cB$ and without loss of generality we assume that $\cB$ is centred at the origin.
Let $\rho_{\p\Omg}\colon\Omg\arr\dR_+$ be the distance function to the boundary of $\Omg$ and let $\rho_{\p\cB}\colon\cB\arr\dR_+$ be the distance function to the boundary of the disk $\cB$. According to, \eg, \cite[Sec. 3]{DZ94} the distance-function $\rho_{\p\Omg}$ is Lipschitz continuous with the Lipschitz constant $=1$, differentiable almost everywhere and
\begin{equation}\label{eq:nablarho}
	|\nb\rho_{\p\Omg}(x)| = 1\qquad \text{for almost all}\,\, x\in\Omg.
\end{equation}
The in-radius of $\Omg$ is defined by
\[
	r_{\rm i} := \max_{x\in\Omg}\rho_{\p\Omg}(x).
\]
It is easy to check by an argument based on the geometric isoperimetric inequality that $r_{\rm i} \le R$ and if $\Omg$ is not congruent to the disk $\cB$ then this inequality is even strict.

For each $t>0$, we define the sub-domains of $\Omg$ and $\cB$ as
\begin{equation}
\begin{aligned}
	\Omg_t &:= \{x\in\Omg\colon\rho_{\p\Omg}(x) > t\},\\
	\cB_t & := \{x\in\cB\colon\rho_{\p\cB}(x) > t\}.
\end{aligned}	
\end{equation}
The lengths of the boundaries of these auxiliary domains satisfy the inequality stated in the next lemma. 
\begin{lem}[{\cite[Prop. A.1]{S01}, \cite{H64}}]\label{lem:levellengths}
For all $t\in(0,r_{\rm i})$,
$|\p\Omg_t| \le L - 2\pi t = |\p\cB_t|$. 
\end{lem}

Our  admissible domains are  those \emph{sub-ordinate} to balls in the sense that the moments of inertia with respect to a fixed center of the level curves of the distance to the boundary are controlled by that for the disk.

\begin{dfn}\label{def:subordinate}
	We say that $\Omg$ is sub-ordinate to $\cB$ if there exists $x_0\in\dR^2$ such that for almost all $t \in (0,r_{\rm i})$ the following inequality holds
	\[
		\int_{\p\cB_t} |x|^2\dd\cH^1(x) = 2\pi(R-t)^3 \ge
		\int_{\p\Omg_t} |x+x_0|^2\dd\cH^1(x),
	\]
	where $\cH^1$ is the one-dimensional Hausdorff measure on the respective curve.
\end{dfn}
The next two propositions give us examples of domains that are sub-ordinate to $\cB$.

\begin{prop}\label{prop:contained}
 If for some $x_0\in\dR^2$ one has  $x_0+\Omg\subset\cB$ then $\Omg$ is sub-ordinate to $\cB$ in the sense of Definition~\ref{def:subordinate}.
\end{prop}
\begin{proof}
First, by Lemma~\ref{lem:levellengths} we have the inequality $|\p\Omg_t| \le |\p\cB_t| =2\pi(R-t)$.
Let $y\in x_0+\Omg_t$ with $t\in(0,r_{\rm i})$ be arbitrary. Hence, we get by a simple geometric argument that $y\in\cB$ and that $\rho_{\p\cB}(y) > t$. Thus, we have the inclusion
$x_0+\Omg_t\subset\cB_t$ for all $t\in (0,r_{\rm i})$.
Hence, for all $x\in \p\Omg_t$ we have $|x+x_0|\le R-t$ and thus the inequality in Definition~\ref{def:subordinate} is satisfied.
\end{proof}
Recall that $\Omg\subset\dR^2$ is said to be \emph{centrally symmetric} 
if it is invariant under the isometric involution $J\colon\dR^2\arr\dR^2$ acting as
$Jx := -x$. \begin{prop}\label{prop:cent-sym}
	If $\Omg\subset\dR^2$ is convex and centrally symmetric then
	it is sub-ordinate to $\cB$ in the sense of Definition~\ref{def:subordinate}. 
\end{prop}
\begin{proof}
	Let $\s = (\s_1,\s_2)\colon[0,\ell]\arr\dR^2$ be the natural parametrization $(|\dot\s(s)| = 1)$ of a piecewise $C^2$-smooth closed curve $\Sg\subset\dR^2$ of length $\ell > 0$.  Assume that the origin is \emph{the centroid} of the curve $\Sg$; \ie
	$\int_\Sg \s(s)\dd s =0$.  Recall that the moment of inertia of
	$\Sg$ with respect to the origin is defined by
	\[
		I_\Sg := \int_\Sg|\s(s)|^2\dd s.
	\]
	Let $\cC\subset\dR^2$ be the circle  of length 
	$\ell > 0$ centred at the origin.
	
	Consider the ordinary differential operator $\sfh\psi :=-\psi''$ with $\dom\sfh:= H^2(\Sg)$ in the Hilbert space $L^2(\Sg)$, which represents the quadratic from $H^1(\Sg)\ni\psi\mapsto \|\psi'\|^2_{L^2(\Sg)}$. The lowest eigenvalue of $\sfh$ is simple, equal to zero and the respective eigenfunction is a constant function. The second eigenvalue of $\sfh$ is equal to $\frac{4\pi^2}{\ell^2}$. 
	Clearly, $\s_1,\s_2\in H^1(\Sg)$ and  applying the min-max principle to the operator $\sfh$
	and using that $\s_1$ and $\s_2$ are both orthogonal to the constant function,
	 we find
	\begin{equation}\label{eq:moment_of_inertia}
		I_\Sg = \int_\Sg (\s_1^2+\s_2^2)\dd s 
		\le 
		\frac{\ell^2}{4\pi^2}\int_\Sg ((\s_1')^2+(\s_2')^2)\dd s
		= \frac{\ell^3}{4\pi^2} = I_\cC.
	\end{equation}
	Since $\Omg$ is centrally symmetric,
	we conclude that $\Omg_t$ is centrally symmetric
	for all $t\in (0,r_{\rm i})$ as well (because $\rho_{\p\Omg}(x)=\rho_{\p\Omg}(-x)$ for all $x\in\Omega$).   Moreover, convexity of $\Omg$ combined with~\cite[Thm. 5.4\,(i)]{DZ94} yields that the distance function $\rho_{\p\Omg}$ is concave in $\Omg$. Hence, $\Omg_t$ is convex and therefore $\p\Omg_t$ is connected.
	It follows by a simple geometric reason that the origin is the centroid of $\p\Omg$ and of the  curves $\p\Omg_t$ for all $t\in(0,r_{\rm i})$.
	Recall also that by \cite[Prop. 6.1]{H64} (see also \cite[Prop. A.1]{S01}) the \vl{(connected)} curve $\p\Omg_t$ is piecewise $C^2$-smooth for almost all $t\in (0,r_{\rm i})$. Hence, combining Lemma~\ref{lem:levellengths} with the inequality~\eqref{eq:moment_of_inertia} we finally obtain that the condition in Definition~\ref{def:subordinate} is satisfied with $x_0 =0$.
\end{proof}
\begin{remark}
	It remains an open question whether there are simply connected $C^2$-smooth domains that are not sub-ordinate to the disk of the same perimeter in the sense of Definition~\ref{def:subordinate}. 
\end{remark}

\begin{remark}
	We remark that the inequality~\eqref{eq:moment_of_inertia} between moments of inertia was first established by Hurwitz~\cite[pp.~396-397]{H02}.  
\end{remark}
\begin{remark}
	In fact, Proposition~\ref{prop:cent-sym} shows slightly more.  A bounded simply connected $C^2$-smooth domain $\Omg\subset\dR^2$ is sub-ordinate
	to the disk $\cB$ of the same perimeter, in the sense of Definition~\ref{def:subordinate},  if the level curves $\p\Omg_t$ are connected and	have the same centroid for almost all $t\in(0,r_{\rm i})$.
\end{remark}
Now we can formulate and prove the main result of this section and of the paper on the isoperimetric inequality for the magnetic Robin Laplacian.
\begin{thm}\label{thm:ipi}
	Let $\Omega$ be a $C^2$-smooth
	bounded simply connected domain sub-ordinate in the sense of Definition~\ref{def:subordinate} to the disk $\cB$ with the same perimeter as $\Omega$.
	Let the set $\cA\subset\dR_-\times\dR_+$ be as in~\eqref{eq:b*-rad}. Let $\lm_1^{\beta,b}(\Omg)$ and $\lm_1^{\beta,b}(\cB)$ be the lowest eigenvalues, respectively, of $\sfH^{\beta,b}_\Omg$ and of  $\sfH_\cB^{\beta,b}$.
	Then for all  $(\beta,b)\in\cA$ the following isoperimetric inequality holds
	\[\lambda_1^{\beta,b}(\Omega)\leq \lambda_1^{\beta,b}(\cB),
	\]
	where the equality occurs if, and only if,
	$\Omg$ is congruent to $\cB$.
\end{thm}
Before giving the proof of the theorem we will formulate its direct corollary, which follows from Theorem~\ref{thm:ipi} combined with Remark~\ref{rem:b*-rad} and Proposition~\ref{prop:estimate_beta}.
\begin{cor}\label{cor:regime-beta-b}
	Let the assumptions be as in Theorem~\ref{thm:ipi}. Let $\beta < 0$ be arbitrary and assume that $0 < b < \min\{R^{-2},4\sqrt{-\beta}R^{-3/2}\}$, where $R>0$ is the radius of the disk $\cB$.
	Then the isoperimetric inequality holds
	\[\lambda_1^{\beta,b}(\Omega)\leq \lambda_1^{\beta,b}(\cB),
	\]
	where the equality occurs if, and only if,
	$\Omg$ is congruent to $\cB$.
\end{cor}
\begin{proof}[Proof of Theorem~\ref{thm:ipi}]
	Without loss of generality we can assume that $\Omg$ is not congruent to the disk $\cB$ and that $\Omg$ is sub-ordinate to the disk $\cB$ in the sense of Definition~\ref{def:subordinate} with $x_0 = 0$.
	In this case we have $r_{\rm i} < R$ where
	as before $r_{\rm i}$ is the in-radius of $\Omg$ and $R$ is the radius of the disk $\cB$.
	
	Let $u_\circ\in H^1(\cB)$ be an eigenfunction
	associated with the ground state for the magnetic Robin Laplacian with the homogeneous magnetic field $b\in\dR_+$ on the disk $\cB$
	and the Robin parameter $\beta$.
	The assumption $(\beta,b)\in\cA$ combined with Proposition~\ref{prop:ef-rad} yields that the eigenfunction $u_\circ$ can be chosen to be a radial function  in the space $C^\infty(\ov\cB;\dR)$, which is positive in $\cB$,  and 	the respective principal eigenvalue, $\lm_1^{\beta,b}(\cB)$, is negative.
	We have the representation $u_\circ(x) = \psi_\circ(\rho_{\p\cB}(x))$ with some 
	$\psi_\circ\in C^\infty([0,R])$, which is positive on $(0,R)$.
	Consider the following test function
	\[
		u_\star(x) := \psi_\circ(\rho_{\p\Omg}(x)),\qquad x\in\Omg.
	\]
	Using Lipschitz continuity of $\rho_{\p\Omg}$ one gets that $u_\star\in H^1(\Omg)$.
	
	Recall that the co-area formula applied in two dimensions, see~\cite[Thm. 4.20]{B19} and~\cite{MSZ02}, to an open set $\cA\subset\dR^2$, a Lipschitz continuous real-valued function $f\colon\cA\arr\dR$, and an integrable function $g\colon\cA\arr\dR$ gives
	\begin{equation}\label{eq:coarea}
	\int_\cA g(x)|\nb f(x)|\,\dd x =
	\int_\dR\int_{f^{-1}(t)} g(x)\,\dd \cH^1(x)\,\dd t,
	\end{equation}
	where $\cH^1$ in the inner integral on the right-hand side is the one-dimensional Hausdorff  measure on the level curve $\{x\in\cA\colon f(x) = t\}$.
	
	In view of~\eqref{eq:nablarho}, we conclude that $|\nb u_\star| = |\psi'_\circ\circ\rho_{\p\Omg}|$ almost everywhere in $\Omg$. Hence, taking that $u_\star$ is real-valued into account, applying the formula~\eqref{eq:coarea} twice to
	$f = \rho_{\p\Omg}$, $\cA = \Omg$, $g = |\nabla u_\star|^2$ in the first term
	and to $f = \rho_{\p\Omg}$, $\cA = \Omg$,
	$g = |x|^2|u_\star|^2$ in the second term below and using again~\eqref{eq:nablarho}, we get
	\begin{equation}\label{eq:kinetic}
	\begin{aligned}
		\|(\nabla& - \ii b {\bf A})u_\star\|^2_{L^2(\Omg;\dC^2)} = \\
		&= \|\nabla u_\star\|^2_{L^2(\Omg;\dC^2)} + \frac{b^2}{4}\int_{\Omg}|x|^2|u_\star(x)|^2\dd x\\
		& =
		\int_0^{r_{\rm i}}|\psi_\circ'(t)|^2\int_{\rho^{-1}_{\p\O\Omg}(t)}\dd\cH^1(x)\dd t\\
		&\qquad\qquad+
		\frac{b^2}{4}
		\int_0^{r_{\rm i}}|\psi_\circ(t)|^2\int_{\rho^{-1}_{\p\O\Omg}(t)}|x|^2\dd\cH^1(x)\dd t
		\\
		& = 
		\int_0^{r_{\rm i}} |\psi_\circ'(t)|^2|\p\Omg_t|\dd t
		 + \frac{b^2}{4}\int_0^{r_{\rm i}}|\psi_\circ(t)|^2\int_{\p\Omg_t}|x|^2 \dd \cH^1(x)\dd t\\
		&<
		\int_0^{R} |\psi_\circ'(t)|^2|\p\cB_t|\dd t + \frac{b^2}{4}\int_0^{R} |\psi_\circ(t)|^2\int_{\p\cB_t}|x|^2 \dd \cH^1(x)\dd t \\
		&= \|(\nabla -\ii b{\bf A})u_\circ\|^2_{L^2(\cB;\dC^2)},\\[0.6ex]
	\end{aligned}
	\end{equation}  
	where in the penultimate step we combined
	that $R > r_{\rm i}$ with the inequality in Lemma~\ref{lem:levellengths} and the inequality in Definition~\ref{def:subordinate} with $x_0 = 0$. 
	
	Using again the co-area formula~\eqref{eq:coarea} and performing the computation analogous to the above we find
	\begin{equation}\label{eq:norm}
			\|u_\star\|^2_{L^2(\Omg)} = \int_0^{r_{\rm i}} |\psi_\circ(t)|^2|\p\Omg_t|\dd t <\int_0^{R} |\psi_\circ(t)|^2|\p\cB_t|\dd t = \|u_\circ\|^2_{L^2(\cB)}.
	\end{equation}			
	Moreover, we obtain that
	\begin{equation}\label{eq:trace}
		\|u_\star|_{\p\Omg}\|^2_{L^2(\p\Omg)} = \|u_\circ|_{\p\cB}\|^2_{L^2(\p\cB)} = L|\psi_\circ(0)|^2.
	\end{equation}
	Combining the min-max principle with~\eqref{eq:kinetic},~\eqref{eq:norm} and~\eqref{eq:trace} and employing the fact that $\lm_1^{\beta,b}(\cB) < 0$ we obtain that
	\begin{align*}
		\lm_1^{\beta,b}(\Omg) &\le \frac{\|(\nabla -\ii b{\bf A})u_\star\|^2_{L^2(\Omg;\dC^2)} +\beta\|u_\star|_{\p\Omg}\|^2_{L^2(\p\Omg)}}{\|u_\star\|^2_{L^2(\Omg)}}\\ 
	&< 
		\frac{\|(\nabla -\ii b{\bf A})u_\circ\|^2_{L^2(\cB;\dC^2)} +\beta\|u_\circ|_{\p\cB}\|^2_{L^2(\p\cB)}}{\|u_\circ\|^2_{L^2(\cB)}} = \lm_1^{\beta,b}(\cB),
	\end{align*}
	where in the last step we used that $u_\circ$ is an eigenfunction of $\sfH_{\cB}^{\beta,b}$ corresponding to its lowest eigenvalue $\lm_1^{\beta,b}(\cB)$.
\end{proof}

\begin{remark}
	Even if there are simply connected $C^2$-smooth planar domains that are not sub-ordinate to the disk of the same perimeter it presents an open question whether the statements of Theorem~\ref{thm:ipi} and of Corollary~\ref{cor:regime-beta-b} hold
	without the subordinacy condition.
\end{remark}

\section{Large coupling asymptotics of $\lm_1^{\beta,b}(\Omg)$ and its connection to the isoperimetric inequality}\label{sec:sa}

We discuss in this section the asymptotics
of $\lm^{\beta,b}_1(\Omg)$ in the limit of large negative Robin parameter ($\beta\arr-\infty$) and its connection with the isoperimetric inequality in Theorem~\ref{thm:ipi}.
Large coupling asymptotics of the lowest Robin eigenvalue in the absence of a magnetic field has been studied by many authors recently~\cite{EMP14, HK17, KP17, P, PP16, LP08}. 
We assume throughout this section that the bounded simply connected domain $\Omg$ is $C^\infty$-smooth and not congruent to the disk. The area of $\Omg$ is denoted by $A$ and the perimeter by $L$. 

It follows from~\cite[Thm. 1]{P} that
\begin{equation}\label{eq:ev-asymp*}
	\lambda_1^{\beta,0}(\Omega)=-\beta^2+\beta\kappa_{\max}(\p\Omega)+\cO(|\beta|^{2/3}),\qquad \beta\arr-\infty,
\end{equation}
where $\kappa_{\max}(\p\Omega)$ is the maximum of the curvature of $\partial\Omega$ and the convention for the sign of the curvature is that the curvature is non-negative for a convex domain.
In the presence of a magnetic field, the ground state energy in the disk $\cB$ of radius $R > 0$ satisfies\footnote{The case  $R\not=1$ can be deduced from the case $R=1$ by a dilation, which yields $\lambda_1^{\beta,b}(\cB)=R^{-2}\lambda_1^{\beta R,bR^2}(\cB_1)$, where $\cB_1$ denotes the unit disk.} (see \cite[Thm. 1.1]{KS}) 
\begin{equation}\label{eq:ev-asymp}
	\lambda_1^{\beta,b}(\cB)=-\beta^2+R^{-1}\beta+R^{-2}e(b,R)+o(1),\qquad\beta\arr-\infty\,,
\end{equation}
where $e(b,R)=-\frac12+\inf_{m\in\Z}\left( m-\frac{bR^2}{2}\right)^{2}=\cO(1)$.

In the general case, the eigenvalue asymptotics $\lambda_1^{\beta,b}(\Omega)$ agrees with \eqref{eq:ev-asymp*} and the contribution of the magnetic field is hidden in the remainder term.
\begin{prop}\label{prop:ev-asymp-b}
For any fixed value of $b\geq 0$, we have, 
\begin{equation}\label{eq:ev-asymp**}
	\lambda_1^{\beta,b}(\Omega)=-\beta^2+\beta\kappa_{\max}(\p\Omega)+\cO(|\beta|^{2/3}),\qquad \beta\arr-\infty.
\end{equation}
\end{prop}
\begin{proof}
By the diamagnetic inequality~\cite[Thm.~7.21]{LL01} and \eqref{eq:ev-asymp*}
\[\lambda_1^{\beta,b}(\Omega)\geq \lambda_1^{\beta,0}(\Omega)= -\beta^2+\beta\kappa_{\max}(\p\Omega)+\cO(|\beta|^{2/3})\,.\]
Consider a normalized and real-valued ground state $u_1^{\beta,0}$ corresponding to the lowest eigenvalue $\lambda_1^{\beta,0}(\Omega)$ of $\sfH_{\Omg}^{\beta,0}$. By the min-max principle and \eqref{eq:ev-asymp*},
\begin{align*}\lambda_1^{\beta,b}(\Omega)&\leq \frq_{\Omega}^{\beta,b}[u_1^{\beta,0}]=
\frq_{\Omega}^{\beta,0}[u_1^{\beta,0}]+b^2\int_\Omg|\Ab|^2|u_1^{\beta,0}|^2\dd x\\
&\leq  \lambda_1^{\beta,0}(\Omega)+\|\Ab\|_\infty^2b^2\\
&\leq  -\beta^2+\beta\kappa_{\max}(\p\Omega)+\cO(|\beta|^{2/3})\,.\qedhere\end{align*} 
\end{proof}
According to~\cite{P2} we have
\begin{equation}\label{eq:curv_ineq}
	\kp_{\rm max}(\p\Omg) > \sqrt{\frac{\pi}{A}} = \kp_{\rm max}\big(\p\cB_{\sqrt{\frac{A}{\pi}}}\big),
\end{equation}
where $\cB_{\sqrt{\frac{A}{\pi}}}$ is the ball of radius $\sqrt{\frac{A}{\pi}}$ and thus having the same area as $\Omg$.
Hence, it follows from the asymptotic expansions~\eqref{eq:ev-asymp} and~\eqref{eq:ev-asymp**} that for a given domain $\Omg$ and $b\ge 0$ there exists a constant $\beta_0(b,\Omega)<0$ such that, for all $\beta\leq \beta_0(b,\Omega)$, 
\[\lambda_1^{\beta,b}(\Omega)< \lambda_1^{\beta,b}\big(\cB_{\sqrt{\frac{A}{\pi}}}\big)\,.\]
Using the geometric isoperimetric inequality $L^2\ge 4\pi A$ we obtain from~\eqref{eq:curv_ineq} 
\[
	\kp_{\rm max}(\p\Omg) > \sqrt{\frac{\pi}{A}} \ge \frac{2\pi}{L} = \kp_{\rm max}(\p\cB_{\frac{L}{2\pi}}),
\]
where $\cB_{\frac{L}{2\pi}}$ is the disk of the radius $\frac{L}{2\pi}$ and thus having the same perimeter as $\Omg$.
Now we can combine the spectral expansions in~\eqref{eq:ev-asymp} and \eqref{eq:ev-asymp**} to deduce that, 
 for $b,\Omega$ fixed, there exists a constant $\beta_1(b,\Omega)<0$ such that, for all $\beta\leq \beta_1(b,\Omega)$,
 \[\lambda_1^{\beta,b}(\Omega)< \lambda_1^{\beta,b}\big(\cB_{\frac{L}{2\pi}}\big)\,,\]
 which is  consistent with the isoperimetric inequality in Theorem~\ref{thm:ipi}.
\begin{remark} 
It is worth to point out that for the above isoperimetric inequalities, which hold for $\beta < 0$ sufficiently large by absolute value, we have not assumed that $\Omg$ is sub-ordinate to the unit disk $\cB$ and that $b$ is moderate.
\end{remark}

\subsection*{Acknowledgement}
AK is partially supported by the Center for Advanced Mathematical Sciences (CAMS, American  University of Beirut).
VL acknowledges the support by the grant No.~21-07129S 
of the Czech Science Foundation (GA\v{C}R)
and thanks Magda Khalile for useful discussions.

\appendix

\section{Closedness and semi-boundedness of the quadratic form $\frq_\Omg^{\beta,b}$}\label{app:clsd_smb}
In this appendix we show that the quadratic form
$\frq_\Omg^{\beta,b}$ in~\eqref{eq:qf} satisfies all the assumptions of the first representation theorem.
\begin{lem}\label{lem:clsd_smb}
	The symmetric densely defined quadratic form $\frq_\Omg^{\beta,b}$ in~\eqref{eq:qf} is closed and semi-bounded.
\end{lem}
\begin{proof}
	Using that $\Ab \in L^\infty(\Omg;\dR^2)$ we find that
	for all $u\in H^1(\Omg)$ one has
	\begin{equation}\label{eq:opposite_bound}
	\begin{aligned}
		\|(\nabla -\ii b\Ab)u\|_{L^2(\Omg;\dC^2)}^2
		&\le 2\|\nabla u \|^2_{L^2(\Omg;\dC^2)} + 2b^2\|\Ab\|_{\infty}^2\|u\|^2_{L^2(\Omg)},\\
		\|(\nabla -\ii b\Ab)u\|_{L^2(\Omg;\dC^2)}^2
		&\ge \frac12\|\nabla u \|^2_{L^2(\Omg;\dC^2)} - b^2\|\Ab\|_{\infty}^2\|u\|^2_{L^2(\Omg)}.
	\end{aligned}	
	\end{equation}
	From the inequalities in~\eqref{eq:opposite_bound} we conclude that the non-negative symmetric densely defined quadratic form $\frq_\Omg^{0,b}$ corresponding to the magnetic Neumann Laplacian on $\Omg$ with the homogeneous magnetic field is closed, because the norm induced by the quadratic form $\frq_\Omg^{0,b}$ is equivalent to the standard norm in the Sobolev space $H^1(\Omg)$.
	
	Recall that according to the diamagnetic inequality~\cite[Thm.~7.21]{LL01}
	\begin{equation}\label{eq:diamag}
	\|\nabla |u|\|_{L^2(\Omg;\dC^2)}^2\le \|(\nabla -\ii b\Ab)u\|_{L^2(\Omg;\dC^2)}^2
	\end{equation}
	for all $u\in H^1(\Omg)$.
	Combining~\eqref{eq:diamag} with the inequality in~\cite[Lem. 2.6]{BEL14} we obtain that for any $\eps >0$ there exists a constant $C(\eps) > 0$ such that
	\begin{equation}\label{eq:s-bd}
		\|u|_{\p\Omg}\|^2_{L^2(\p\Omg)} \le \eps\|(\nabla -\ii b\Ab)u\|^2_{L^2(\Omg;\dC^2)} + C(\eps)\|u\|^2_{L^2(\Omg)},\quad \text{for all}\, u\in H^1(\Omg).
	\end{equation}
	From the above inequality we deduce that the quadratic form $H^1(\Omg)\ni u\mapsto \beta\|u|_{\p\Omg}\|^2_{L^2(\p\Omg)}$, $\beta\in\dR_-$, is form bounded with respect to the quadratic form $\frq^{0,b}_\Omg$ with the form bound $< 1$.
	Hence, by \cite[Thm. VI.1.33]{K} the quadratic form $\frq^{\beta,b}_\Omg$ is closed and semi-bounded.
\end{proof}

\section{Continuity of the ground state energy}\label{sec:cont}

In this appendix we present a standard proof   that the ground state energy of the magnetic Robin Laplacian $\sfH_{\Omg}^{\beta,b}$  depends continuously on the intensity of the magnetic
field $b$ and the Robin parameter $\beta$.

Recall that $\Omg\subset\dR^2$ is a bounded simply connected $C^2$-smooth domain.
For definiteness we use the convention $\|u\|^2_{H^1(\Omg)} := \|\nabla u\|^2_{L^2(\Omg;\dC^2)} + \|u\|^2_{L^2(\Omg)}$ for the standard norm in the Sobolev space $H^1(\Omg)$. Recall also that by the trace theorem~\cite[Thm. 3.38]{McL} there exists a constant $c  >0$ such that $\|u|_{\p\Omg}\|^2_{L^2(\p\Omg)}\le c\|u\|^2_{H^1(\Omg)}$ for any $u\in H^1(\Omg)$.

Let $\beta_1\le0$ and $b_1\ge 0$ be fixed and $\beta_2\le0$ and $b_2\ge 0$ be such that $|\beta_1-\beta_2|,|b_1-b_2|\le1$.
It follows from the second inequality in~\eqref{eq:opposite_bound} combined with~\cite[Lem. 2.6]{BEL14} that there exists $\gamma\in\dR$ such that $\sfH^{\beta_2,b_2}_\Omg \ge \gamma$ for
any $\beta_2\le 0$ and $b_2\ge 0$ satisfying $|\beta_1-\beta_2|,|b_1-b_2|\le1$

For any $u \in H^1(\Omg)$, we get the following estimate
\[
\begin{aligned}
&\big|\frq^{\beta_1,b_1}_\Omg[u] - \frq^{\beta_2,b_2}_\Omg[u]\big|\\  
&\quad\!=\!\bigg|2 (b_1-b_2)\Im (\nabla u,{\bf A} u)_{L^2(\Omg;\dC^2)}
+ (b_1^2-b_2^2)(|{\bf A}|^2 u,u)_{L^2(\Omg)}  + (\beta_1-\beta_2)\|u|_{\p\Omg}\|^2_{L^2(\p\Omg)} 
\bigg|\\
&\quad\!\le\!
2|b_1-b_2|\|\nabla u\|_{L^2(\Omg;\dC^2)}
\|{\bf A} u\|_{L^2(\Omg;\dC^2)}
\!+\! |b_1^2-b_2^2|\|{\bf A}\|_\infty^2\|u\|^2_{L^2(\Omg)}\!+\!  
|\beta_1-\beta_2|\|u|_{\p\Omg}\|^2_{L^2(\p\Omg)}\\
&\quad\!\le\!
|b_1\!-\!b_2|\big[\|\nabla u\|_{L^2(\Omg;\dC^2)}^2
\!+\!\|{\bf A}\|^2_\infty \|u\|^2_{L^2(\Omg)}\big]
\!+\!|b_1^2\!-\!b_2^2|\|{\bf A}\|_\infty^2\|u\|^2_{L^2(\Omg)}\! +\! c|\beta_1\!-\!\beta_2|\|u\|^2_{H^1(\Omg)}\\
&\quad \!\le\! \max\{|b_1-b_2|,|b_1-b_2|
\|{\bf A}\|_\infty^2, |b_1^2-b_2^2|\|{\bf A}\|_\infty^2,c|\beta_1-\beta_2|\}\|u\|^2_{H^1(\Omg)},
\end{aligned}	
\] 
where we used the trace theorem in the penultimate step.
Since the standard $H^1$-norm is equivalent to the norm $u \mapsto \frq_{\Omg}^{\beta_1,b_1}[u] +(1-\gamma)\|u\|^2_{L^2(\Omg)}$ induced by the quadratic form $\frq_\Omg^{\beta_1,b_1}$ we conclude from the above estimate with the aid of~\cite[Thm. VI.3.6]{Kato} that the operator $\sfH^{\beta_2,b_2}_\Omg$ converges in the norm resolvent sense to the operator $\sfH^{\beta_1,b_1}_\Omg$ as $(\beta_2,b_2)\arr(\beta_1,b_1)$. Note also that the family of operators $\sfH^{\beta_2,b_2}_\Omg$ is uniformly lower semibounded.  Hence, it follows from the spectral convergence result~\cite[Satz 9.24\,(ii)]{W00} that $\lm^{\beta_2,b_2}_1(\Omg)\arr \lm^{\beta_1,b_1}_1(\Omg)$ as $(\beta_2,b_2)\arr(\beta_1,b_1)$ and thus the lowest eigenvalue $\lm_1^{\beta,b}(\Omg)$ of $\sfH^{\beta,b}_\Omg$ is a continuous function of the parameters $\beta,b\in(-\infty,0]\times [0,\infty)$.

 \section{The Neumann magnetic ground state energy}\label{app:C}

Consider the Neumann eigenvalue $\lambda_1^{0,b}(\Omega)$ introduced in \eqref{eq:ev} with   an associated normalized eigenfunction $u_b\colon\Omg\arr\dC$.
Let us assume that $\lambda_1^{0,b}(\Omega)=0$.  The diamagnetic inequality, 
\[ \int_\Omega\big|\nabla|u_b|\big|^2\dd x\leq \int_\Omega|(\nabla-\ii b\Ab)u_b|^2 \dd x=0,\]
yields that  $|u_b|= |\Omega|^{-1/2}$,  since the domain $\Omega$ is connected.
Furthermore,  $u_b$ satisfies
\begin{equation}\label{eq:grad-m-ub} 
(\nabla-\ii b \Ab)u_b=0\,.
\end{equation}
Taking the inner product with $\overline{u}_b$, we infer from \eqref{eq:grad-m-ub},
\[ \overline{u}_b\nabla u_b=\frac{\ii b}{|\Omega|}\Ab\,.\]
Taking the ${\rm curl}$ in \eqref{eq:grad-m-ub},  we get,
\begin{equation}\label{eq:grad-m-ub*}
 \ii b\,{\rm curl}(\Ab u_b)={\rm curl}(\nabla u_b)=0 \,.\end{equation}
Finally, we notice that (see \eqref{eq:mp})
\[\overline{u}_b\,{\rm curl}(\Ab u_b)=|u_b|^2\,{\rm curl}\Ab+\Ab^\bot\cdot(\overline{u}_b\nabla u_b)=\frac{1}{|\Omega|}+\frac{\ii b}{|\Omega|}\Ab^\bot\cdot\Ab=\frac{1}{|\Omega|}\]
where $\Ab^\bot=\frac12(x_1,x_2)$. Consequently,  we get from \eqref{eq:grad-m-ub*} that
 $b=0$.

\end{document}